\documentclass[10pt,oneside,reqno]{amsart}
\usepackage[cp1251]{inputenc}
\usepackage[english]{babel}
\usepackage{amsmath}
\usepackage{amssymb}
\usepackage{amsfonts}

\usepackage{hyperref}
\newtheorem{Th}
{Theorem}
\newtheorem{Lem}[Th]{Lemma}

\newtheorem{Remark}[Th]{Remark}

\begin{document}
\sloppy


\begin{center}

\textbf{Alternative algebras admitting derivations with invertible values and invertible derivations.}
\footnote{The authors were supported 
by RFBR 12-01-31016, 12-01-33031,
by RF President Grant council for support of young scientists and leading scientific schools (project MK-330.2013.1),
and FAPESP (Grant 2011/51132-9).}

\medskip

\medskip
\textbf{Ivan Kaygorodov$^{a,b}$, Yury Popov$^{b,c}$}

\medskip

$^a$ Instituto de Matem\'{a}tica e Estat\'{i}stica, Universidade de S\~{a}o Paulo, Brasil,\\
$^b$ Sobolev Institute of Mathematics, Novosibirsk, Russia,\\
$^c$ Novosibirsk State University, Russia.
\sloppy
\sloppy

\medskip

\tableofcontents 

\medskip

\section{Introduction.}
\end{center}

\medskip
 
The notion of derivation with invertible values as a derivation of a ring with unity that 
takes only multiplicatively invertible or zero values appeared in \cite{Berg}.
Bergen, Herstein and Lanski determined the structure of associative rings that admit derivations with invertible values.
Later, the results of this paper were generalized in \cite{Giam}--\cite{Komatsu}.

Another interesting type of derivations are invertible derivations. 
The definition of an invertible derivation as an invertible mapping first arose in \cite{Jac}, 
where the nilpotency of a Lie algebra admitting an invertible derivation was proved.
The research on that topic was then continued in \cite{Baj,Moens}. 

Nowadays, a great interest is shown in the study of nearly associative algebras and superalgebras with derivations.
For example, works \cite{Popov,Popov2} determine the structure of differentiably simple alternative and Jordan algebras, 
and papers \cite{Filll}--\cite{kay_okh} 
give the description of generalizations of derivations of simple and semisimple alternative, Jordan and structurable (super)algebras.
Nevertheless, the problem of specification of algebras from classical non-associative varieties (such as alternative, Jordan, structurable, etc.), 
admitting derivations with invertible values and invertible derivations, remains unconsidered.
The present work is to make up this gap.

\medskip

\section{Basic definitions and identities.}

\medskip

We are using standard notation:
$$(x,y,z):=(xy)z-x(yz)\mbox{ --- the associator of elements }x, y, z,$$ 
$$\left[x,y\right]:=xy-yx\mbox{ --- the commutator of elements }x, y,$$
$$x\circ y:=xy+yx\mbox{ --- the Jordan product of elements }x, y.$$

An algebra $A$ is called \textit{alternative} 
(see \cite{kolca} for more information on alternative algebras), if $A$ satisfies the following identities:
\begin{eqnarray*} (x,x,y)=0, (x,y,y)=0. \end{eqnarray*}
It's easy to check that in any alternative algebra the associator is a skew-symmetric function of its arguments, 
and the \textit{flexible identity} $x(yx) = (xy)x$ holds.
It's also well known  \cite[p.35]{kolca} that every alternative algebra satisfies the \textit{middle Moufang identity}: $(xy)(zx) = x(yz)x$.
A commutative algebra $J$ is called \textit{Jordan} if it satisfies the \textit{Jordan identity}:
\begin{eqnarray*} (x^2,y,x) = 0. \end{eqnarray*}
It is widely known that if $A$ is an alternative algebra, 
then vector space $A$ with new multiplication $a \circ b$ is a Jordan algebra which we will denote by $A^{(+)}$.

The \textit{nucleus} of an algebra $A$ is the set
\begin{eqnarray*}N(A) = \left\{n \in A |\ (n,A,A)=(A,n,A)=(A,A,n)=(0)\right\},\end{eqnarray*}
the \textit{commutative center} of $A$ is the set 
\begin{eqnarray*}K(A)=\left\{k \in A|\ \left[k,A\right]=\left[A,k\right]=(0)\right\}, \end{eqnarray*}
 and the \textit{center} of $A$ is  $Z(A)=N(A)\cap K(A)$.

Derivation $d$ is called \textit{inner} if it lies in the smallest subspace of the space of all linear operators on $A$ 
containing all right and left multiplications by elements of $A$ and closed under commutation. Otherwise $d$ is called \textit{outer}.

In studying the structure of alternative algebras, one class is of great importance: Cayley--Dickson algebras. 
The definition and properties of Cayley--Dickson algebras and the Cayley--Dickson process can be found, for instance, in \cite{kolca}. 
It's known that every Cayley--Dickson algebra $C$ over field $F$ is $8$--dimensional, non-associative, alternative, simple and has an unit element.
Also, $C$ is \textit{quadratic} over $F$, that is, for every $x \in C$ the following relation holds:
\begin{eqnarray}
\label{odin}
 x^2-t(x)x+n(x)=0, 
 \end{eqnarray}
where $t(x),\ n(x)\in F,\ t(x)$ is a $F$-linear mapping, and $n(x)$ is a strictly nondegenerate quadratic form satisfying 
$n(xy)=n(x)n(y)$ for all $x, y \in C.$ 
A Cayley--Dickson algebra is 
also equipped with a symmetric bilinear nondegenerate form $f(x,y) = n(x+y) - n(x) - n(y).$ 
For a subset $M \subseteq C,$ by $M^{\bot}$ we will denote the orthogonal complement to $M$ with respect to $f.$

A Cayley--Dickson algebra containing zero divisors is called \textit{split}. It's known \cite[p.43]{kolca} that element $x$ of a 
split Cayley--Dickson algebra is invertible if and only if $n(x) \neq 0$.

It's also known \cite[p.46]{kolca} that every split Cayley--Dickson algebra over field $F$ is isomorphic to a 
\textit{Cayley--Dickson matrix algebra} $C(F)$, comprising matrices of the form $a = \begin{pmatrix} \alpha & u \\ v & \beta \end{pmatrix}$, where $\alpha, \beta \in F$, $u, v \in F^3$.

Addition and scalar multiplication of elements of the algebra $C(F)$ will then correspond to the usual addition and scalar multiplication of matrices. 
However, multiplication of elements of the algebra $C(F)$ will correspond to the following matrix multiplication:

\begin{center}
$\begin{pmatrix} \alpha & u \\ v & \beta \end{pmatrix} \cdot \begin{pmatrix} \gamma & t \\ w & \delta \end{pmatrix} = \begin{pmatrix} \alpha\beta + (u,w) & \alpha t + \delta u - v\times w \\ \gamma v + \beta w + u \times t & \beta \gamma + (v,t) \end{pmatrix}$,
\end{center}
where for vectors $x=(x_1,x_2,x_3), y=(y_1,y_2,y_3) \in F^3,$ by 
$$(x,y)=x_1 y_1 + x_2 y_2 + x_3 y_3$$ we denote their dot product, and by
$$x \times y = (x_2 y_3 - x_3 y_2, x_3 y_1 - x_1 y_3, x_1 y_2 - x_2 y_1)$$
their cross product.
Under given representation $t(a) = \alpha + \beta$, $n(a)= \alpha \beta - (u,v)$.

In the case when $char (F) \neq 2$, $C$ can be obtained from $F$ by applying the Cayley--Dickson process thrice to $F$ with the identical involution 
and parameters $\alpha, \beta, \gamma \in F$. We will not go into the full details here and will only provide the formula that defines 
multiplication in algebra $B = A + vA$ obtained by the Cayley--Dickson process from algebra $A$ with involution \ $\bar{}$ \ :
\begin{eqnarray*} (a_1+vb_1)(a_2+vb_2) = (a_1a_2+\gamma b_2\overline{b_1})+v(\overline{a_1}b_2+a_2b_1), \end{eqnarray*}
where $a_i, b_i \in A$, $v^2 = \gamma \in F$.

We will also need the following statement, which describes simple alternative non-associative algebras.

\begin{Th}
\label{Th1.1}
Let $A$ be a simple non-associative alternative algebra.
Then the center of the algebra $A$ is a field and $A$ is a Cayley--Dickson algebra over its center.
\end{Th}

\medskip

\section{Alternative algebras with derivations with invertible values.} 

\medskip
Let $A$ be an algebra with unit element $1$ over field $F$. We will denote by $U$ the set of invertible elements of $A$.
In this section we will only consider \textit{derivations with invertible values}, 
by which we understand such non--zero derivations $d$ that for every $x \in A,$ $d(x) \in U$ or $d(x)=0$  holds.

In 1983, Bergen, Herstein and Lanski initated the study whose purpose is to relate the structure of a ring to the special behavior 
of one of its derivations. Namely, in their article \cite{Berg} they described associative rings admitting derivations with invertible values. 
They proved that such ring must be either a division ring, or the ring of $2 \times 2$ matrices over a division ring, 
or a factor of a polynomial ring over a division ring of characteristic $2$. They also characterized those division rings 
such that a $2 \times 2$ matrix ring over them has an inner derivation with invertible values.
Further, associative rings with derivations with invertible values (and also their generalizations) were discussed in variety of works 
(see, for instance, \cite{Giam}--\cite{Komatsu}).
So, in \cite{Giam}, semiprime associative rings with involution, allowing a derivation with invertible values 
on the set of symmetric elements, were given an examination. 
In work  \cite{Car} Bergen and Carini determined the associative rings admitting a derivation with invertible values on some non--central Lie ideal. 
Also, in papers \cite{Chang} and \cite{Hongan} the structure of associative rings 
that admit $\alpha$-derivations with invertible values 
and their 
natural generalizations --- $(\sigma,\tau)$-derivations with invertible values --- was described. 
And in paper \cite{Komatsu} Komatsu and Nakajima described associative rings that allow generalized derivations with invertible values. 

The purpose of this section is to generalize the results of Bergen, Herstein and Lanski to the alternative case.

\medskip
In this part, $A$ is an alternative algebra with unit element $1$ and derivation with invertible values $d$.
The following lemmas were proved in \cite{Berg} for associative algebras and can be easily generalized to the alternative case 
with minor differences, but in order to ensure the complteness of the narration we shall provide their proofs.

\medskip

\begin{Lem}
\label{Lem1.2} 
If $d(x)=0$, then either $x=0$, or $x$ is invertible.
\end{Lem}
\textbf{Proof.} Let's notice \cite[p.204]{kolca}, that in every alternative algebra the following identity holds:
\begin{eqnarray}
 \label{dva}
(a^{-1},a,b)=0.
\end{eqnarray}  
It's then easy to see that in an arbitrary alternative algebra the product of two invertible elements is also invertible. 
Using identity (\ref{dva}), for invertible $a$ and $b$ we find
$$(b^{-1}a^{-1})(ab)= a^{-1}((ab)b^{-1}) - (a^{-1}(ab))b^{-1} +(b^{-1}a^{-1})(ab)=$$
$$-(a^{-1},ab,b^{-1}) + (b^{-1}a^{-1})(ab)= -(b^{-1},a^{-1},ab) + (b^{-1}a^{-1})(ab)=
b^{-1}(a^{-1}(ab))=1.$$
Assume that $x\neq0$. Since $d \neq 0$, there exists $y \in A$ such that $d(y) \in U$. 
Hence $d(yx)=d(y)x \in U$ and $d(y)^{-1}d(yx)=x$. 
In view of $d(y)$ and $d(yx)$ being invertible, $x$ is also invertible. 
The lemma is proved.

\medskip

Now we shall study the ideal structure of $A$:
\begin{Lem}
\label{Lem1.3}
$a)$ If $L \neq 0$ is a one--sided ideal in $A$ then $d(L) \neq (0)$.\\
$b)$ If $I$ is a proper one--sided ideal of $A$, then $I$ is both minimal and maximal.\\
$c)$ If $I$ is a proper ideal of $A$ then $I^2=(0)$.\\
$d)$ If $char(A) \neq 2$ then $A$ is simple.
\end{Lem}

\medskip

\textbf{Proof.} 
$a)$ Since the statement is obvious when $L=A$, we should only consider the case when $L \neq A$. If $0 \neq a \in L$, then, 
by lemma \ref{Lem1.2}, 
$d(a) \neq 0$, since $a$ is not invertible.\\
$b)$ It suffices to show that every proper one--sided ideal in $A$ is maximal. 
Let $I \subset J$ be a proper one--sided ideal in $A$. It's easy to see that $d(I) \cap I=(0)$ and $I \oplus d(I)$ is also an one--sided ideal in $A$. 
By lemma \ref{Lem1.3}$(a)$, $d(I) \neq (0)$, hence $d(I)$ contains invertible elements, in consequence of which  $I \oplus d(I)=A$. For arbitrary $j \in J $ we have $j=a+d(b)$, $a,b \in I$. Consequently, $d(b)=j-a \in J \cap d(I) = (0)$; thus $j=a \in I$.\\
$c)$ If $I \neq A$ is an ideal of $A$, then
$$d(I^2) \subset d(I)I+Id(I) \subset I,$$ 
consequently, 
by lemma \ref{Lem1.3}$(a)$, $I^2=(0)$, since the product of two ideals in an alternative algebra is also an ideal \cite[p.115]{kolca} 
and $I$ does not contain any invertible elements.\\
$d)$ Let $2A \neq 0$ and $I \neq (0)$. Then, by lemma \ref{Lem1.3}$(a)$, there exists $b \in I$ such that $d(b) \in U$. Since $b^2=0$, 
\begin{eqnarray*}0=d^2(b^2)=d^2(b)b + 2d(b)^2 + bd^2(b), \end{eqnarray*} 
and consequently $2d(b)^2 \in I$.
Now, since $d(b)$ is invertible, $d(b)^2$ is also invertible and $2d(b)^2 = 0$, therefore $2=0$. 
We have obtained a contradiction which proves the lemma.

\medskip

By $Der(A)$ we will denote the set of all derivations of algebra $A$. Let us fix some subset $D \subseteq Der(A)$.
The ideal $I$ is called a \textit{$D$--ideal}, if for all $\partial \in D$,  $x \in I$ we have $\partial(x) \in I$.
Algebra $A$ is called \textit{$D$--simple} if $A^2 \neq 0$ and $A$ contains no proper $D$--ideals 
(for more detailed information on $D$--simple algebras see \cite{Popov,Popov2} and their references).

\medskip

As an immediate consequence of lemma \ref{Lem1.3}$(a)$ we have

\begin{Lem}
If alternative algebra $A$ admits a derivation with invertible values $d$, then $A$ is $d$--simple.
\end{Lem}
\medskip

Now, if $char (A) \neq 2$ we can apply lemma \ref{Lem1.3}$(d)$ and theorem \ref{Th1.1} and 
conclude that $A$ is either an associative or Cayley--Dickson algebra over its center. 
We will now consider the non--simple non--associative case, which is examined in

\begin{Lem}
\label{Lem2.3}
If $A$ is not simple and not associative,
then $A = C[x]/(x^2)$, where $C$ is a Cayley--Dickson algebra over its center $Z(C)$, 
$C$ is a division algebra, $char(C) = 2$, $d(C)=0$, $d(x)=1+ax$ for some $a \in Z(C)$, and $d$ is an outer derivation.
\end{Lem}
 
\medskip

\textbf{Proof.} Combining lemma \ref{Lem1.3}$(b)$ and $(d)$, we have $char(A)=2$, $I^2=(0)$ for any proper ideal $I$ in $A$ and all
proper one--sided ideals in $A$ are both minimal and maximal. Consequently, we can easily deduce that $A$ contains a unique 
(left, right, two--sided) ideal $M$ and $M^2=0$. 
Therefore, as in the proof of lemma \ref{Lem1.3}$(b)$, we have $A=M \oplus d(M)$, 
particularly, for any $a \in A$ there exist $m,n \in M$ such that $d(a)=m+d(n)$. 
Hence $m=d(a-n) \in M \cap d(A)=(0)$ and so, denoting $C=\ker (d)$, we have $A=C + M$. 
By lemma \ref{Lem1.2}, $C$ is a division algebra, therefore $A=C \oplus M$. 
We define linear mappings $\lambda:M \rightarrow C $ and $\mu:M \rightarrow M$ by $d(m)=\lambda(m)+\mu(m)$ for any $m \in M$. 
It's easy to notice that for any $a \in C, b \in M$ the following holds: 
$$a\mu(b)+a\lambda(b)=ad(b)=d(ab)=\mu(ab)+\lambda(ab),$$ 
where $a\mu(b), \mu(ab) \in M$ and consequently $a\lambda(b)=\lambda(ab) \in \lambda(M)$;
similarly $\lambda(ba)=\lambda(b)a \in \lambda(M)$. This implies that $\lambda(M)$ is an ideal in $C$. 
Since $C$ is simple and $\lambda(M) \neq (0)$ we derive that $C$ is isomorphic to $M$ as a left $C$--module. Putting $x=\lambda^{-1}(1)$, we have $A=C \oplus Cx$. Using the fact that $\lambda$ is a module isomorphism, it's easy to see that $[x,C]=0$. Considering the identity
\begin{eqnarray*}3(k,x,y)=3(y,k,x)=3(x,y,k)=[xy,k]-x[y,k]-[x,k]y=0,\end{eqnarray*}
satisfied for any $k \in K(B)$, $x,y \in B$ in arbitrary alternative algebra $B$ \cite[p.136]{kolca}, and taking into account the structure of $A$ we deduce that $x \in Z(A)$. Therefore we have $A \cong C[x]/(x^2)$. Now nonassociativity of $A$ and theorem \ref{Th1.1} imply that $C$ is a Cayley--Dickson algebra over its center $Z(C)$. 
We can write $\mu(x)=ax$ for some $ a \in C$.
Now, since $x \in Z(A)$ and $char (A) = 2$, for arbitrary $c \in C$ we have:
\begin{eqnarray*}0=d(cx+xc)=c(1+ax)+(1+ax)c=cax+axc=(ca+ac)x. \end{eqnarray*}
Since $C$ is a division algebra, we obtain $ca+ac=0$, thus $a \in Z(C)$.

Finally, since every ideal of $A$ is invariant under the action of any inner derivation, $x \in M$, and $d(x)\notin M$, it is clear that $d$ is not inner. 
The lemma is proved.
\medskip 

\begin{Th}
Let $A$ be an alternative algebra with unit element $1$, admitting derivation with invertible values $d$. 
Then$:$
\\$1)$ $A$ is an associative algebra and one of the following conditions holds:
\\$a)$ $A$ is a division algebra $D$;
\\$b)$ $A$ is a $2 \times 2$ matrix algebra $M_2 (D)$ over division algebra $D$;
\\$c)$ $A$ is a factor--algebra of polynomial algebra $D[x]/(x^2)$ over division algebra $D$; 
furthermore, $char (D)=2,\ d(D)=0$ and $d(x)=1+ax$ for some $a$ in the center of $D,$ and $d$ is an outer derivation;\\
$2)$ $A$ is a non--associative alternative algebra and one of the following conditions holds:
\\$a)$ $A$ is a Cayley--Dickson algebra over its center $Z(A)$;
\\$b)$ $A$ is a factor--algebra of polynomial algebra $C[x]/(x^2)$ over a Cayley--Dickson division algebra; 
furthermore, $char(C)=2,\ d(C)=0$ and $d(x)=1+ax$ for some $a$ in the center of $C,$ and $d$ is an outer derivation.
\end{Th}

\medskip

\textbf{Proof.} 
The associative case follows from \cite{Berg}, and the non-associative case follows from theorem \ref{Th1.1}, lemmas \ref{Lem1.3} and \ref{Lem2.3}.

\medskip
Now, to complete the characterization of alternative algebras allowing derivations with invertible values we only 
have to describe split Cayley--Dickson algebras with derivations with invertible values, which is done in the following.

\begin{Lem}
\label{Lem2.2} 
An algebra $C$, which is a split Cayley--Dickson algebra over its center $Z$, admits a derivation with invertible values $d$ 
if and only if one of the following conditions holds:\\
$I)$ $C$ is obtained by means of the Cayley--Dickson process from its associative division subalgebra $B$: 
$C=B+vB, v^2=\gamma \in Z$, $\gamma \neq 0$, where $B = \ker (d)$ and  $dim_Z B = 4$. 
Furthermore, in this case an arbitrary derivation with invertible values $d$ is of the form $d(a+vb)=v(bu)$, 
where $a, b \in B$ and $u \in B$ is a fixed element with $t(u)=0$.\\
$II)$ $C$ can be represented as a direct sum: $C=B+xB$, where  $t(x)=0$, $B = \ker (d)$, $B$ is a subfield of $C$, 
$B=B^{\bot}$ and $dim_Z B = 4$. 
Furthermore, in this case an arbitrary derivation with invertible values $d$ is of the form $d(a+xb)=b$, where $a, b \in B$.
\end{Lem}
\medskip

\textbf{Proof.} It's generally known (see, for example, \cite{Eld}) that every derivation of $C$ is inner. 
It's easy to see then that $Z \subseteq ker (d)$ and $d$ is a $Z$--linear mapping. Therefore we will consider $C$ as a $Z$--algebra.
Suppose that $C$ allows a derivation with invertible values $d$. 
Take a subspace $V \subset C$ such that $dim_Z V = 4$ and $V$ does not contain invertible elements. For example (taking into account that $C \cong C(F)$ --- the Cayley--Dickson matrix algebra over $F$), we can take
\begin{center}
$V=\left\{\begin{pmatrix} \alpha & u \\ 0 & 0\end{pmatrix}\ | \alpha \in Z, u \in Z^3\ \right\}.$
\end{center}
From lemma \ref{Lem1.2} it follows that $dim_Z d(V) = 4$ and $V \cap d(V) = (0)$, hence $C=V \oplus d(V)$. 
 In particular, for any $x \in C$ there exist $u, v \in V$ such that $d(x) = u + d(v)$. 
 Consequently, $u=d(x-v) \in V \cap d(A) = (0)$, and, denoting $B=\ker (d)$, we have $C=B + V$. 
 By lemma \ref{Lem1.2}, $B$ is a division algebra, thus $C=B \oplus V$ and $dim_Z B = 4$. 
 Combining the facts that $B$ is simple, $Z(C) \subseteq Z(B)$ and applying theorem \ref{Th1.1} we have that $B$ is an associative subalgebra in $C$.
 By \cite[p.39]{kolca}, in $C$ the following relation is valid: 
\begin{eqnarray}
\label{tri}
a \circ b -t(a)b -t(b)a-f(a,b)=0.
\end{eqnarray}
Putting $b=d(a)$, we obtain
\begin{eqnarray}
\label{chet}
a \circ d(a) -t(a)d(a) -t(d(a))a-f(a,d(a))=0.
\end{eqnarray}
Applying $d$ on (\ref{odin}), we have
\begin{eqnarray}
\label{pyat}
a \circ d(a) -t(a)d(a)=0.
\end{eqnarray}
Subtracting (\ref{chet}) from (\ref{pyat}), we obtain $t(d(a))a+f(a,d(a))=0$.
If $a$ and 1 are linearly independent over $Z$, then we have
\begin{eqnarray}
\label{shest}
f(a,d(a))=0.
\end{eqnarray}
In the case when $a \in Z$, then $a \in \ker (d)$ and relation (\ref{shest}) is then obvious.
Linearizing (\ref{shest}), we obtain $f(a,d(b))+f(d(a),b)=0$. 
Consequently, since $B=\ker(d)$, for arbitrary $a \in C$ we have $f(d(a),B)=-f(a,d(B))=0,$ 
and so $d(C) \subseteq B^{\bot}$. We now have to study two cases:

$(I)$ If the restriction of the form $f$ on $B$ is nondegenerate, then, by \cite[p.32, Th.1]{kolca}, 
$C$ can be obtained from $B$ by means of the Cayley--Dickson process, that is, $C = B + vB$, $v^2 = \gamma \neq 0$, $B^{\bot}=vB$. 
Particularly, $d(v)=vu$ for some $u \in B$, and therefore for arbitrary $a, b \in B$ we have $d(a+vb)=d(v)b=(vu)b=v(bu)$.

By \cite[p.26]{kolca}, for any $x, y, w \in C$ we have
\begin{eqnarray*}
n(x)f(y,w)=f(xy,xw).
\end{eqnarray*}
And for $x=v$, $y=1$, $w=u$, using (\ref{shest}), we obtain $$0=f(v,vu)=n(v)t(u).$$ 
Since $v^2=\gamma \in Z$, $\gamma \neq 0$, then $n(v)\neq0$ and $t(u)=0$.

$(II)$ Now, let the restriction of the form $f$ on $B$ be degenerate. Hence there exists $0\neq b \in B$ such that $f(b,B)=0$. Therefore $0=f(b,b)=2n(b)$. 
Since $b$ is invertible, $n(b)\neq 0$ and we must have $char (C) = 2$.
By \cite[p.26]{kolca}, in $C$ the following relation holds:
\begin{eqnarray}
\label{sem}
f(x,z)f(y,w)=f(xy,zw)+f(xw,yz).
\end{eqnarray}
Putting into  (\ref{sem}) $x=b$, $z=a$, $y=b^{-1}c$, $w=1$, where $a,c \in B$, we have 
$$0=f(b,a)f(b^{-1}c,1)=f(c,a)+f(b,ab^{-1}c),$$ 
and so by arbitrariness of $a,c$ we conclude that $f(B,B)=0$, that is, $B \subseteq B^{\bot}$.

Now we shall show that opposite inclusion also takes place:
Suppose for a moment that there exists $x \in B^{\bot}$, $x \notin B$. 
By (\ref{dva}) and skew--symmetry of the associator, $dim_Z xB = 4$ and $A = B \oplus xB$.
Using (\ref{sem}), we have $$f(a,xc)=f(a\cdot1,xc) = -f(ac,x) + f(a,x)f(1,c)=0$$
for any $a,c \in B$. Consequently, $xB \subset B^{\bot}$ and $C=B^{\bot}$, which contradicts the nondegeneracy of the form $f$.
We put $x=d^{-1}(1)$. It's obvious that $x \notin B$ and $C = B \oplus xB$.
Relation (\ref{shest}) implies that $0=f(x,1)=t(x)$.
Now we only have to prove that $B$ is a field. 
By the definition of $f$ and the fact that $f(B,B) = 0$, for any $a,c$ we have 
$0 = f(a,c)=n(a+c)-n(a)-n(c)$, which means that $n$ is a ring homomorphism from $B$ to $Z$. 
In view of $B$ being simple, together with $n(1)=1$, $ker (n) = 0$, and we conclude that $B$ is a subfield of $Z$.

{\bf Conversely}, 
suppose that condition $(I)$ holds, that is $C$ is obtained  from $B$ by means of the Cayley--Dickson process.
Let $0 \neq u$ be an element of $B$ such that $t(u)=0$. Consider the mapping $d:a+vb \mapsto v(bu)$, where $a, b \in B$.
We are to show that $d$ is a derivation. Indeed, for any $a_1, b_1, a_2, b_2 \in B$ we have 
$$d(a_1+vb_1)(a_2+vb_2)+(a_1+vb_1)d(a_2+vb_2)=$$
$$\gamma(b_2(u+\overline{u})\overline{b_1}) + v((a_2b_1+\overline{a_1}b_2)u)=\gamma(b_2t(u)\overline{b_1}) + v((a_2b_1+\overline{a_1}b_2)u)=$$
$$v((a_2b_1+\overline{a_1}b_2)u)=d((a_1a_2+\gamma b_2\overline{b_1})+v(\overline{a_1}b_2+a_2b_1))=d((a_1+vb_1)(a_2+vb_2)).$$
Also, $$n(d(a+vb))=n(v(bu))=n(v)n(b)n(u)= -\gamma n(b)n(u)\neq 0$$ 
if $b \neq 0$, since $B$ is a division algebra. 
Hence $d(a+vb)$ is invertible for any $a \in B$, $0 \neq b \in B$, so $d$ takes invertible values.

Now assume that condition $(II)$ holds.
Consider the mapping $d:a+xb \mapsto b$. We are to show that $d$ is a derivation with invertible values.
Since we have $B=B^{\bot}$, then for any $a \in B,$ $t(a)=0$ holds. Combining (\ref{tri}) and $char (C) = 2$, we obtain 
\begin{eqnarray}
\label{vosem}
\left[x,a\right] = x\circ a = t(a)x + t(x)a + f(a,x) = f(a,x) \in Z,
\end{eqnarray}
particularly, $d(\left[x,a\right]) = 0$.
Substituting $x$ in (\ref{odin}), we deduce that $x^2 \in Z$.
Using (\ref{vosem}), it's easy to check that for $a, c \in B$ the following identity holds:
\begin{eqnarray}
\label{devyat}
(a,c,x) = af(c,x)+ f(a,x)c + f(x,ac),
\end{eqnarray}
and consequently, $d((a,c,x)) = 0$. Now we will prove that $d$ is a derivation. 
For arbitrary $a, b, c, h \in B$ we have 
$$d((ax+b)(cx+h)) = d((ax)(cx) + (ax)h + b(cx) + bh) = d((ax)(cx)) + d((ax)h) + d(b(cx)).$$ 
Consider the last two summands: 
$$d((ax)h) = d((xa)h)=d(x(ah))=ah, \ d(b(cx)) = d((bc)x) = bc.$$ 
On the other hand, 
$$d(ax+b)(cx+h) + (ax+b)d(cx+h) = a(cx+h) + (ax+b)c = a(cx) + ah + (ax)c + bc.$$ 
Therefore we need to show that $d((ax)(cx)) = a(cx) + (ax)c$. 
Transforming the corresponding expressions, we have: 
$$a(cx) + (ax)c = (ac)x + (a,c,x) + a(xc) + (a,x,c) = (ac)x + a(cx + f(c,x)) = (a,c,x) + af(c,x).$$ 
Using the middle Moufang identity, we obtain 
$$d((ax)(cx)) = d((xa + f(a,x))cx) = d((xa)(cx)) + f(a,x)d(cx) = d(x(ac)x) + f(a,x)c =$$
$$d(x(x(ac) + f(x,ac))) + f(a,x)c  =d(x^2(ac)) + f(x,ac) + f(a,x)c = f(x,ac) + f(a,x)c,$$ 
since $x^2 \in Z$ and $d(x^2(ac))=n(x)d(ac)=0$. 
Equating the expressions, we will arrive at the relation (\ref{devyat}), which, as was shown earlier, holds identically. 
Therefore $d$ is a derivation of $C$. Since $d$ takes values in $B,$ which is a field, it's obvious that $d$ is a derivation with invertible values. 
The lemma is proved.
\medskip

\textbf{Example.} In work \cite{Gagola} an example of a split Cayley--Dickson algebra $C$ which has a subfield of dimension 4 was provided. Let's consider an imperfect field $F$ of characteristic 2 and elements $\alpha, \beta \in F$ such that $\alpha, \beta, \alpha\beta$ are linearly independent over $F^2$. Then subalgebra $B$ of matrix Cayley--Dickson algebra $C(F)$, generated by elements
\begin{center}
$\begin{pmatrix} 0 & (\alpha,0,0) \\ (1,0,0) & 0 \end{pmatrix}$, $\begin{pmatrix} 0 & (0,\beta,0) \\ (0,1,0) & 0 \end{pmatrix}$,
\end{center}
is a subfield of $C$, and $dim_F B = 4$.

\medskip

\section{A characterization of nilpotent alternative algebras by invertible Leibniz-derivations.} 

\medskip 

In 1955, Jacobson \cite{Jac} proved that a Lie algebra over a
field of characteristic zero admitting a non--singular (invertible)
derivation is nilpotent. The problem of whether the inverse of this
statement is correct remained open until work \cite{Dix}, where
an example of nilpotent Lie algebra, whose derivations are
nilpotent (and hence, singular), was constructed. Such types of
Lie algebras are called characteristically nilpotent Lie algebras.

The study of derivations of Lie algebras leads to the appearance of the notion of their
natural generalization --- a pre-derivation of a Lie algebra, which is a derivation of a Lie triple system induced by that algebra. 
In \cite{Baj} it was proved that Jacobson's result is
also true in terms of pre-derivations. 
Several examples of nilpotent Lie
algebras whose pre-derivations are nilpotent were presented in
\cite{Baj}, \cite{Burd}.

In paper \cite{Moens} a generalization of derivations and
pre-derivations of Lie algebras is defined as a Leibniz-derivation of
order $k$. 
Moens proved that a Lie algebra over a field of characterisic zero is nilpotent if and only if
it admits an invertible Leibniz-derivation. 
After that, 
Fialowski, Khudoyberdiyev and Omirov \cite{Fialc12} showed
that with the definition of Leibniz-derivations from \cite{Moens}
the similar result for non-Lie Leibniz algebras is not true.
Namely, they gave an example of a non--nilpotent Leibniz algebra which
admits an invertible Leibniz-derivation. In order to extend the
results of paper \cite{Moens} for Leibniz algebras they introduced a
definition of Leibniz--derivations of Leibniz algebras which agrees
with the case of Leibniz-derivations of Lie algebras,  and proved
that a Leibniz algebra is nilpotent if and only if it admits an
invertible Leibniz-derivation. 
It should be noted that there exist non-nilpotent Filippov ($n$-Lie) algebras with invertible derivations (see \cite{Williams}). 
Also, in \cite{kay_pre} a generalization of pre-derivations of associative algebras was considered.

The main purpose of this section is to prove the analogue of 
Moens's theorem for alternative algebras.
Throughout the section all spaces of algebras are assumed finite-dimensional over a field of characteristic zero.

\medskip 

{\bf Definition.}
A Leibniz-derivation (by Moens) of order $n$ for an algebra $A$ is an endomorphism $\phi$
of that algebra satisfying the identity
$$\phi((\ldots (x_1x_2) \ldots)x_n)=\sum\limits_{i=1}^n(\ldots (\ldots (x_1x_2) \ldots \phi(x_i)) \ldots)x_n.$$

\medskip 

\begin{Th}
\label{Th2.1}
An alternative algebra over a field of characteristic zero is nilpotent if and only if it has 
an invertible Leibniz-derivation.
\end{Th}

\medskip 

{\bf Proof.}
Let $A$ be a finite--dimensional alternative algebra with an invertible Leibniz-derivation $\phi$ of order $n$ 
and $\beta(A)$ be the nilpotent radical of $A$ (it's also widely known that in the finite-dimensional case it coincides with $rad(A),$ the solvable radical of $A$).
Using \cite{kolca}, we can establish that $A/\beta(A)$ can be represented as finite sum of its minimal ideals, 
where each of them is either a full matrix algebra over some division ring or a Cayley--Dickson algebra over its center. 
Therefore, algebra $A/\beta(A)$ possesses unit element $1$. We will regard $A$ as a direct sum: 
$A=A^s +\beta(A),$ where $A^s$ is a semisimple alternative algebra isomorphic to $A/\beta(A)$ (Wedderburn--Malcev decomposition).
Using the idea of the proof from \cite{Moens} we shall prove that $\phi(\beta(A))\subseteq \beta(A)$. 
We will remark that in the case when $\phi$ is a derivation it was proved 
for all algebras with locally nilpotent radical in \cite{Slinko1972}.

{\bf Step 1.}
We define on vector space $A$ the structure of $n$-ary algebra $A_n$ with multiplication
$$[a_1, a_2,\ldots, a_n]_n=(\ldots (a_1a_2) \ldots )a_n.$$
Hence $\phi$ is a derivation of $n$-ary algebra $A_n$.
We shall show that solvable radicals $rad(A_n)$ and $rad(A)$ of algebras $A_n$ and $A$ coincide.
It's clear that $rad(A) \subseteq rad(A_n)$. 
Consider the natural projection $\pi : A \rightarrow A^s.$
It's easy to see that $\pi(rad(A_n))$ is a solvable ideal in $A^s$: 
applying $\pi$ to the both sides of relation
$$[A, \ldots, rad(A_n), \ldots, A]_n \subseteq rad(A_n),$$
and using the fact that $A^s$ has a unit, we have
$$\pi (rad(A_n)) A^s + A^s \pi (rad(A_n)) \subseteq \pi (rad(A_n)). $$
Consequently, since $A^s$ is semisimple, we have $\pi(rad(A_n))=0.$

{\bf Step 2.} 
We will now show that $\phi(\beta(A))\subseteq \beta(A)$.
Let $\beta(A)=\tau=\tau_1=rad(A_n)$ and $\tau_{t+1}=[\tau_t, \tau_t, \ldots, \tau_t]_n.$ Then we have
$$\tau=\tau_1 \supsetneq \tau_2 \supsetneq \ldots \supsetneq \tau_p \supsetneq \tau_{p+1} =0.$$
Since the product of two ideals in an alternative algebra is also an ideal, then $\tau_t$ is an ideal in $A_n$ for any $t$.

We need to show that $\phi^i(\tau_t) \subseteq \tau$ holds for any $i$.
We use induction on $t$. The induction base is trivial for $t=p+1.$
Now let's suppose that $\phi^i (\tau_{t+1}) \subseteq \tau$ for arbitrary $i$.
We need to prove that the inclusion $\phi^i (\tau_{t}) \subseteq \tau$ holds for any $i$.

The set $\tau+\phi(\tau_t)$ is a solvable ideal of $A_n$, since
$$[A,\ldots, A, \tau+\phi(\tau_t), A, \ldots, A]_n \subseteq \tau +\tau_t+\phi(\tau_t) =\tau+\phi(\tau_t)$$
and
$$[\tau+\phi(\tau_t), \ldots, \tau+\phi(\tau_t)]_n 
\subseteq \tau +[\phi(\tau_t), \ldots, \phi(\tau_t)]_n \subseteq \tau + \phi^{n}(\tau_{t+1}) \subseteq \tau.$$

Now we are to show that $\tau+\phi^k(\tau_t)$ is a solvable ideal of $A_n$ for any $k$.
Suppose that $\phi^i(\tau_t) \subseteq \tau$ for each $1\leqslant i<k$.

Using the induction hypothesis, we have
$$[A,\ldots, A, \tau+\phi^k(\tau_t), A, \ldots, A]_n \subseteq $$
$$\tau + \phi([A,\ldots, A, \phi^{k-1}(\tau_t), A, \ldots, A]_n)+ \sum [ A,\ldots, \phi(A), \ldots, A, \phi^{k-1}(\tau_t), A, \ldots, A]_n) \subseteq \ldots \subseteq $$

$$\tau +
\sum_{ a_0 + \ldots +a_{n-1} =k, a_i \geq 0} \phi^{a_0} ([ \phi^{a_1}(A), \ldots,  \phi^{a_{l-1}}(A), \tau_t, \phi^{a_l}(A), \ldots, \phi^{a_{n-1}}(A)]_n) 
+ \phi^k([A,\ldots, A, \tau_t, A, \ldots, A]_n)
\subseteq $$
$$\tau + \sum
_{i=0}^{k-1} \phi^i(\tau_t) + \phi^k(\tau_t) =\tau +\tau_t+\phi^k(\tau_t) =\tau+\phi^k(\tau_t)$$
and
$$[\tau+\phi^k(\tau_t), \ldots, \tau+\phi^k(\tau_t)]_n \subseteq \tau +[\phi^k(\tau_t), \ldots, \phi^k(\tau_t)]_n \subseteq  
\tau + \phi^{kn}(\tau_{t+1}) \subseteq \tau.$$
Therefore,
$\phi^i(\tau_t) \subseteq \tau$ and $\phi(\tau) \subseteq \tau.$

{\bf Step 3.} 
Considering the fact that $\phi$ is an invertible mapping, that is, it has a trivial kernel,
we conclude that $\phi(A/\beta(A))=A/\beta(A)$, which contradicts the unitality of $A/\beta(A)$, since $\phi(1)=n\phi(1)$ and $\dim(\phi(A/\beta(A)))<\dim(A/\beta(A)).$
This contradiction implies that $A=\beta(A)$, that is, $A$ is nilpotent.

{\bf Step 4.} 
The converse also takes place: in order to see that nilpotent alternative algebra $A$ with nilpotency index $s$
has an invertible Leibniz--derivation of order $n=[\frac{s}{2}]+1$ it suffices to consider the sum of vector spaces $A=W+A^n$
and linear mapping $\phi,$ defined this way:

$$\phi(x)=\left\{\begin{array}{ll}
x,  \mbox{ if } x \in W,\\
nx, \mbox{ if } x \in A^n.
\end{array} \right.$$

It is easy to see that $\phi$ is a Leibniz-derivation of $A$ of order $n$.
The theorem is proved.

\medskip 

\medskip 

Further, it's easy to check that

\medskip

\begin{Remark}
\label{Rem2.2}
Over a field of positive characteristic there exist a nilpotent alternative algebra possessing only singular derivations.
\end{Remark}

\medskip

{\bf Proof.}
Non--associative alternative nilpotent algebras of dimension not greater than $7$ were classified in \cite{Badalov1984}.
Using this classification, over a field of characteristic $3$ we define a $7$--dimensional algebra $A$ 
with basis $\{e_1,e_2,e_3,u_1,u_2,v,w\}$ by this mulptiplication table (all other products are zero): 
$$e_1^2=u_1, e_2^2= u_2, 
e_2e_3=e_3e_2=-v,e_3e_1=u_2, e_1e_3=u_2,$$ 
$$e_1u_1=u_1e_1=v,e_2u_2=u_2e_2=w, e_1v=ve_1=u_1^2=w.$$
It's easy to notice that $A^2=\langle u_1, u_2, v,w\rangle, A^3=\langle v,w\rangle, A^4=\langle w\rangle.$
Then for any derivation $D$ of $A$ we have 
$$D(v)=D(e_1e_1e_1) \in A^4 \mbox{ and }D(w)=D(e_1e_1e_1e_1) \in A^4,$$
then there exist $x \neq 0$ such that $D(x)=0.$

\medskip

We also note that the theorem only holds for characteristic zero:
\begin{Remark}
\label{Rem2.3}
For an arbitrary alternative algebra over a field of positive characteristic $p$
the identity map is a Leibniz-derivation of order $p+1$.
\end{Remark}

\medskip 

\begin{Remark}
\label{Rem2.4}
The free alternative algebra with $n$ generators admits an invertible derivation, but it is not nilpotent. 
\end{Remark}

\medskip 

{\bf Proof.} It suffices to consider a derivation which acts identically on the generators of the algebra.

\medskip 

\begin{Remark}
Theorem \ref{Th2.1}, Remarks \ref{Rem2.3} and   \ref{Rem2.4} take place in case of associative algebras too.
\end{Remark}

\medskip

\begin{Remark}
Following the article \cite{Slinko1972} and using methods provided in the proof of Theorem \ref{Th2.1}
we can show that a finite-dimensional Jordan algebra admitting an invertible derivation is nilpotent.
\end{Remark}

\medskip

\section{Alternative and Jordan algebras with \texorpdfstring{$QDer=End.$}{Qder=End.}} 

\medskip

Following Leger and Luks \cite{LL}, 
we call an additive mapping $f$ a quasiderivation if there exists a linear map $Q$ such that
$Q(xy)=f(x)y+xf(y).$
Leger and Luks described all finite--dimensional Lie algebras in which an arbitrary endomorphism is a quasiderivation.
They found out that such an algebra is either an abelian Lie algebra, a two--dimensional solvable Lie algebra or a three--dimensional simple Lie algebra.
Later on, this result was generalized for Lie superalgebras in \cite{ZZ}.
Also, quasiderivations and generalized derivations were studied in \cite{kg11}--\cite{kay_lie2} and in other works.

\medskip 

Here we shall describe all finite-dimensional Jordan and alternative algebras over arbitrary fields of characteristic not equal to 2,
in which every endomorphism is a quasiderivation.

\begin{Th}
Let $J$ be a finite--dimensional Jordan algebra over a field of characteristic not equal to 2 and such that $QDer(J)=End(J).$ 
Then either $J$ is a field or $J$ has zero multiplication.
\end{Th}

{\bf Proof.} Let $f$ be an endomorphism of algebra $J$. 
Hence there exists $Q_f$ such that $Q_f(xy)=f(x)y+xf(y).$
Suppose that $0 \neq x \in Ann(J) = \{a \in J | a\cdot J = 0\}.$ 
Then for arbitrary $y \in J$ we have $f(x)y=Q_f(xy)-xf(y)=Q_f(xy)=0$. 
Since $f(x)$ can be any element of algebra $J$, then either $Ann(J)=0$ or $J$ has zero multiplication.

Suppose for a moment that there exists $x$ such that $x^2=0$. 
Then for any $f \in End(J)$ we have $2f(x)x=Q(x^2)=0$ 
and by the above argument $J$ has zero multiplication.
It's now easy to notice, if we suppose that $J$ has a nonzero nilpotent radical, then considering that $J$ is a power--associative 
algebra we obtain that there exists a nilpotent element of index 2. 
Therefore we conclude that $J$ has a trivial nilpotent radical and consequently $J$ can be represented as a direct sum of simple Jordan algebras.
Notice (see, for example, \cite{shestakov}) that every direct summand is $f$--invariant, 
thus from now on $J$ will be regarded as a simple unital Jordan algebra.
Description of quasiderivations of simple finite-dimensional Jordan algebras \cite{shestakov} 
implies that $f$ can be represented as a sum of a scalar mapping and a derivation of algebra $J$.
We notice that if we denote by $R_a$ the operator of right multiplication by element $a \in J$,
then $R_J + Der(J)$ is a subspace in $End(J)$.
Furthermore, $R_J \cap Der(J) =0,$ since for every unital algebra $1R_a=a$ and $dim(R_J)=dim(J)$.
But $dim(Der(J)) \neq dim(J)^2-1$, if $dim(J)\neq 1.$
That is, we can conclude that $J$ is a field.
The theorem is now proved.

\begin{Th}\label{Thfin}
Let $A$ be a finite-dimensional alternative algebra over a field of characteristic not equal to 2 and such that $QDer(A)=End(A).$ 
Then either  $A$ is a field or $A$ has zero multiplication.
\end{Th}

{\bf Proof.} It suffices to notice that if $f$ is a quasiderivation of algebra $A$ then $f$ is a quasiderivation of algebra $A^{(+)}$.
Furthermore, if $A$ is an alternative algebra then $A^{(+)}$ is a Jordan algebra. 
Consequently, if $A$ satisfies the condition $QDer(A)=End(A),$ then $A^{(+)}$ satisfies this condition too, therefore it is either a field or a zero multiplication algebra.
Therefore we can conclude that $A$ is either a field or an anticommutative algebra.
It's generally known that an alternative anticommutative algebra is nilpotent, that is, it has nontrivial annihilator $Ann(A)$. 
So, if we have $0 \neq x \in Ann(A)$ then for any $y \in A$ we have $f(x)y=0$ and by the above argument $A$ has zero multiplication.
The theorem is proved.

\medskip

\begin{Remark}
Theorem \ref{Thfin} for associative algebras can be proved as consequence of the result of Leger and Luks \cite{LL}.
It suffices to notice that a quasiderivation $f$ of algebra $A$ is also a quasiderivation of algebra $A^{(-)},$
and to consider associative commutative algebras and associative algebras $A$ such that $A^{(-)}$ 
is either a two--dimensional solvable Lie algebra or a three--dimensional simple Lie algebra.
\end{Remark}

\medskip

{\bf Acknowledgements.}
The authors are grateful to Prof. Ivan Shestakov  (IME-USP, Brazil)
and Prof. Alexandre Pozhidaev (Sobolev Inst. of Math., Russia) for interest and constructive comments,
and Mark Gannon (IME-USP, Brazil) for the translation.

\end{document}